\documentclass{article}[12pt]
\usepackage{amssymb,amsfonts}

\begin{document}

\begin{center}
\Large A note on the convergence of multivariate formal \\ power series solutions of meromorphic Pfaffian systems
\end{center}
\medskip
\centerline{\large I.\,V.\,Goryuchkina, R.\,R.\,Gontsov}
\bigskip
\begin{abstract}
Here we present some complements to theorems of Gerard and Sibuya, on the convergence of multivariate formal power series solutions
of nonlinear meromorphic Pfaffian systems. Their the most known results concern completely integrable systems with nondegenerate
linear parts, whereas we consider some cases of non-integrability and degeneracy.
\end{abstract}

\section{Introduction}

Consider a Pfaffian system
\begin{equation}\label{Psyst}
dy={\bf f}_1(x,y)\frac{dx_1}{x_1^{p_1}}+\ldots+{\bf f}_m(x,y)\frac{dx_m}{x_m^{p_m}}, \qquad x=(x_1,\ldots,x_m),\quad
y=(y_1,\ldots,y_n)^{\top},
\end{equation}
where ${\bf f}_i:\,({\mathbb C}^{m+n},0)\rightarrow({\mathbb C}^n,0)$ are germs of holomorphic maps and $p_i\geqslant1$
are integers. Equivalently, this system has the form
$$
\Theta:=dy-\omega=0,
$$
where $\omega$ is a (${\mathbb C}^n$-valued) differential 1-form meromorphic in a neighbourhood $D$ of
$0\in{\mathbb C}^{m+n}$, with the polar locus
$$
\Sigma=\{(x,y)\in D\mid h(x,y):=x_1\ldots x_m=0\}.
$$
In the case $p_1=\ldots=p_m=1$ we have the Pfaffian system (\ref{Psyst}) of {\it Fuchs type}, in this case $\omega$ is
a {\it logarithmic} 1-form in $D$, that is, $h\,\omega$ and $h\,d\omega$ are holomorphic in $D$.

Written in a PDEs form, (\ref{Psyst}) becomes
\begin{equation}\label{PDE}
x_1^{p_1}\,\frac{\partial y}{\partial x_1}={\bf f}_1(x,y), \quad\ldots,\quad
x_m^{p_m}\,\frac{\partial y}{\partial x_m}={\bf f}_m(x,y).
\end{equation}
We study a question of the convergence of a formal power series solution
\begin{equation}\label{series}
\varphi=\sum_{|k|>0}{\bf c}_kx^k\in{\mathbb C}[[x]]^n \quad ({\bf c}_k\in{\mathbb C}^n)
\end{equation}
of such a system. Here, as usual in a multivariate case,
$$
k=(k_1,\ldots,k_m)\in{\mathbb Z}_+^m, \quad |k|=k_1+\ldots+k_m, \quad x^k=x_1^{k_1}\ldots x_m^{k_m}.
$$
A basic work on this subject (as well as on the analytic and asymptotic properties of such series) is that by R.\,Gerard and
Y.\,Sibuya \cite{GS}, in which Pfaffian systems are assumed to be {\it completely integrable} on $D\setminus\Sigma$.
This means that for any $(x^0,y^0)\in D\setminus\Sigma$, there exists a unique solution $y=y(x)$ of (\ref{Psyst}) such that
$y(x^0)=y^0$. Due to the Frobenius theorem, the complete integrability of (\ref{Psyst}) is equivalent to the relation
$$
d\Theta=\Omega\wedge\Theta,
$$
for some matrix differential 1-form $\Omega$ holomorphic in $D\setminus\Sigma$ (see \cite[Ch.1, Th. 5.1]{IKSY} or \cite{War}). However,
in \cite{GS} there are obtained some assertions concerning the convergence of (\ref{series}) that don't apply to the complete integrability
of (\ref{Psyst}) and thus hold for any Pfaffian system of the form (\ref{Psyst}). These are the following two theorems.
\smallskip

{\bf Theorem A.} {\it Any formal power series solution $(\ref{series})$ of a Fuchsian system
\begin{equation}\label{Fuchs}
x_1\,\frac{\partial y}{\partial x_1}={\bf f}_1(x,y), \quad\ldots,\quad
x_m\,\frac{\partial y}{\partial x_m}={\bf f}_m(x,y)
\end{equation}
converges near $0\in{\mathbb C}^m$.}
\smallskip

{\bf Theorem B.} {\it If $p_j=1$ for some j $($that is, the system $(\ref{PDE})$ is Fuchsian along the component $\{x_j=0\}$ of its polar locus$)$
and the corresponding Jacobi matrix $\partial{\bf f}_j/\partial y(0,0)$ does not have non-negative integer eigenvalues then
any formal power series solution $(\ref{series})$ of $(\ref{PDE})$ converges near $0\in{\mathbb C}^m$ $($in the case of the complete
integrability of such a system, there holds the existence and uniqueness of the solution$)$.}
\smallskip

We will show that the sufficient condition of convergence from Theorem B can be weakened for a fixed formal solution $\varphi$ in the
following way (see the proof in Section 3).
\smallskip

{\bf Theorem 1.} {\it Let $p_1=1$ in the system $(\ref{PDE})$ and let $\varphi$ be its formal power series solution. If the matrix
$$
A=\frac{\partial{\bf f}_1}{\partial y}(x,\varphi)|_{x_1=0}\in{\rm Mat}(n,{\mathbb C}[[x_2,\ldots,x_m]])
$$
is such that $\det(A-jI)\not\equiv0$ for any non-negative integer $j$, then $\varphi$ converges near $0\in{\mathbb C}^m$.}
\smallskip

In the following theorem from \cite{GS} the complete integrability of the Pfaffian system is required.
\smallskip

{\bf Theorem C.} {\it In the non-Fuchsian case with all $p_i>1$, if there are $j\ne l$ such that the Jacobi matrices
$\partial{\bf f}_j/\partial y(0,0)$, $\partial{\bf f}_l/\partial y(0,0)$ are non-degenerate then there is a unique formal power series
solution $(\ref{series})$ of $(\ref{PDE})$ and it converges near $0\in{\mathbb C}^m$.}
\smallskip

In a more recent paper by Sibuya \cite{Si} (where $m=2$) the assumption of complete integrability in the above Theorem C is omitted
and there is proved that any formal power series solution $(\ref{series})$ of $(\ref{PDE})$ converges near $0\in{\mathbb C}^m$,
under the rest assumptions of the theorem.

In Section 2 we discuss the condition of complete integrability in more details and study the convergence of formal power series
solutions of $(\ref{PDE})$ in the non-integrable case (Theorems 2, 3). In Section 3 we give the proof of Theorem 1, and in Section 4
we complement the above Gerard--Sibuya Theorem C for the non-Fuchsian case by some sufficient condition of the convergence of
(\ref{series}) satisfying (\ref{PDE}) with the zero Jacobi matrices $\partial{\bf f}_i/\partial y(0,0)$.


\section{The relations of complete integrability}

If the system (\ref{PDE}) is completely integrable, for any its solution $y(x)$ the equality of the second partial derivatives
$\partial^2y/\partial x_i\partial x_j$ and $\partial^2y/\partial x_j\partial x_i$ implies
$$
\frac1{x_i^{p_i}}\frac{\partial}{\partial x_j}{\bf f}_i(x,y(x))=
\frac1{x_j^{p_j}}\frac{\partial}{\partial x_i}{\bf f}_j(x,y(x)),
$$
whence
$$
\frac1{x_i^{p_i}}\frac{\partial {\bf f}_i}{\partial x_j}(x,y(x))-\frac1{x_j^{p_j}}\frac{\partial {\bf f}_j}{\partial x_i}(x,y(x))=
\frac1{x_i^{p_i}x_j^{p_j}}\frac{\partial {\bf f}_j}{\partial y}(x,y(x)){\bf f}_i(x,y(x))-
\frac1{x_i^{p_i}x_j^{p_j}}\frac{\partial {\bf f}_i}{\partial y}(x,y(x)){\bf f}_j(x,y(x)).
$$
Since in the case of complete integrability for any $(x^0,y^0)\in D\setminus\Sigma$ there exists a unique solution $y=y(x)$ of
(\ref{PDE}) such that $y(x^0)=y^0$, one has
\begin{equation}\label{integr}
\frac1{x_i^{p_i}}\frac{\partial {\bf f}_i}{\partial x_j}-\frac1{x_j^{p_j}}\frac{\partial {\bf f}_j}{\partial x_i}\equiv
\frac1{x_i^{p_i}x_j^{p_j}}\Bigl(\frac{\partial {\bf f}_j}{\partial y}{\bf f}_i-\frac{\partial {\bf f}_i}{\partial y}{\bf f}_j\Bigr),
\quad i,j=1,\ldots,m,
\end{equation}
in $D\setminus\Sigma$. Or, equivalently, all the (vector) functions
$$
{\bf F}_{ij}(x,y)=x_j^{p_j}\frac{\partial {\bf f}_i}{\partial x_j}-x_i^{p_i}\frac{\partial {\bf f}_j}{\partial x_i}+
\frac{\partial {\bf f}_i}{\partial y}{\bf f}_j-\frac{\partial {\bf f}_j}{\partial y}{\bf f}_i,
\quad i,j=1,\ldots,m,
$$
are equal to zero identically. Conversely, let us show that if all the functions ${\bf F}_{ij}(x,y)$, $i,j=1,\ldots,m$, equal zero
identically then the system (\ref{PDE}) is completely integrable. Indeed, in this case (\ref{integr}) holds and for
$\Theta=dy-\sum_{i=1}^m{\bf f}_i(x,y)dx_i/x_i^{p_i}$ we have
\begin{eqnarray*}
d\Theta&=&\sum_{i=1}^m\sum_{j=1}^m\frac1{x_i^{p_i}}\frac{\partial {\bf f}_i}{\partial x_j}dx_i\wedge dx_j+
\Bigl(\sum_{i=1}^m\frac1{x_i^{p_i}}\frac{\partial {\bf f}_i}{\partial y}dx_i\Bigr)\wedge dy=
\sum_{i<j}\Bigl(\frac1{x_i^{p_i}}\frac{\partial {\bf f}_i}{\partial x_j}-\frac1{x_j^{p_j}}\frac{\partial {\bf f}_j}{\partial x_i}\Bigr)dx_i\wedge dx_j+\\
& & +\Omega\wedge dy=\sum_{i<j}\frac1{x_i^{p_i}x_j^{p_j}}\Bigl(\frac{\partial {\bf f}_j}{\partial y}{\bf f}_i-
\frac{\partial {\bf f}_i}{\partial y}{\bf f}_j\Bigr)dx_i\wedge dx_j+\Omega\wedge dy=\Omega\wedge\Theta,
\end{eqnarray*}
where
$$
\Omega=\sum_{i=1}^m\frac1{x_i^{p_i}}\frac{\partial {\bf f}_i}{\partial y}dx_i
$$
is a matrix differential 1-form holomorphic in $D\setminus\Sigma$. Hence, the Frobenius integrability condition is fulfilled.

Thus we see that the complete integrability of the system (\ref{PDE}) is described by at most $m(m-1)/2$ vector
($nm(m-1)/2$ scalar) relations (since ${\bf F}_{ij}=-{\bf F}_{ji}$). Further we use the two results by A.\,Ploski \cite{Pl1}
(a version with the detailed proof is \cite{Pl2}) following from his sharpened version of Artin's Approximation Theorem \cite{Art}.
These are:

1) {\it if $f(x_1,\ldots,x_m,y)$ is a non-zero germ of a holomorphic function $({\mathbb C}^{m+1},0)\rightarrow({\mathbb C},0)$
and $\varphi\in{\mathbb C}[[x_1,\ldots,x_m]]$ is a formal power series without constant term such that $f(x_1,\ldots,x_m,\varphi)=0$,
then $\varphi$ converges near $0\in{\mathbb C}^m$;}

2) {\it if ${\bf f}(x_1,\ldots,x_m,{\bf y})$ is a germ of a holomorphic map $({\mathbb C}^{m+n},0)\rightarrow({\mathbb C}^n,0)$,
$\varphi\in{\mathbb C}[[x_1,\ldots,x_m]]^n$ is a $($vector$)$ formal power series without constant term such that ${\bf f}(x_1,\ldots,x_m,\varphi)=0$,
and
$$
\det\frac{\partial{\bf f}}{\partial{\bf y}}(x_1,\ldots,x_m,\varphi)\ne0,
$$
then $\varphi$ converges near $0\in{\mathbb C}^m$.}
\medskip

{\bf Theorem 2.} {\it In the case of the scalar unknown $y$ $($i.e., $n=1)$, if the system $(\ref{PDE})$ is not completely integrable then
any its formal power series solution $\varphi$ converges near $0\in{\mathbb C}^m$.}
\medskip

{\bf Proof.} Non-integrability implies that at least one of the (scalar in this case) functions ${\bf F}_{ij}(x,y)$ is not identically zero
(say, ${\bf F}_{12}$), while ${\bf F}_{12}(x,\varphi)=0$. Thus the assertion follows from the first result by Ploski.
\medskip

{\bf Theorem 3.} {\it Let the system $(\ref{PDE})$ be non-completely integrable and have a formal power series solution
$\varphi\in{\mathbb C}[[x_1,\ldots,x_m]]^n$. If among all the vectors ${\bf F}_{ij}(x,y)$ there are $n$ non-zero components $g_1,\ldots,g_n$
such that
$$
\det\frac{\partial g_i}{\partial y_j}(x_1,\ldots,x_m,\varphi)\ne0,
$$
then $\varphi$ converges near $0\in{\mathbb C}^m$.}
\medskip

{\bf Proof.} Follows from the second result by Ploski.

\section{Proof of Theorem 1}

Clearly, the formal solution $\varphi$ is represented in the form
$$
\varphi=\sum_{j=0}^{\infty}{\bf c}_j(x_2,\ldots,x_m)\,x_1^j, \quad {\bf c}_j\in{\mathbb C}[[x_2,\ldots,x_m]]^n.
$$
Let us prove the convergence of all ${\bf c}_j$ in some common polydisc $\Delta\subset{\mathbb C}^{m-1}$ centered at the origin.

From the equality
\begin{equation}\label{p1}
x_1\,\frac{\partial\varphi}{\partial x_1}={\bf f}_1(x,\varphi)
\end{equation}
it follows that
$$
{\bf f}_1(0,x_2,\ldots,x_m,{\bf c}_0)=0.
$$
Since by the theorem assumption
$$
\det A=\det\frac{\partial{\bf f}_1}{\partial y}(0,x_2,\ldots,x_m,{\bf c}_0)\ne0,
$$
the series ${\bf c}_0$ converges near $0\in{\mathbb C}^{m-1}$ due to the second result by Ploski and thus the matrix $A$ is holomorphic
in some polydisc $\Delta\subset{\mathbb C}^{m-1}$ centered at the origin.

Further we represent $\varphi={\bf c}_0+\psi:={\bf c}_0+\sum_{j=1}^{\infty}{\bf c}_j\,x_1^j$ in the equality (\ref{p1}) and obtain
$$
x_1\,\frac{\partial}{\partial x_1}\Bigl(\sum_{j=1}^{\infty}{\bf c}_j\,x_1^j\Bigr)={\bf f}_1(x,{\bf c}_0+\psi),
$$
that is,
\begin{equation}\label{cj}
\sum_{j=1}^{\infty}j{\bf c}_j\,x_1^j={\bf f}_1(x,{\bf c}_0)+\sum_{j=1}^{\infty}\frac{\partial{\bf f}_1}{\partial y}(x,{\bf c}_0){\bf c}_j\,x_1^j+O(\psi^2),
\end{equation}
whence
$$
(A-I){\bf c}_1=-\frac{\partial{\bf f}_1}{\partial x_1}(0,x_2,\ldots,x_m,{\bf c}_0).
$$
More generally, if we already have that ${\bf c}_0,\ldots,{\bf c}_{j-1}$ are holomorphic in $\Delta$ then (\ref{cj}) implies
$$
(A-jI){\bf c}_j={\bf h}_j(x_2,\ldots,x_m,{\bf c}_0,\ldots,{\bf c}_{j-1}),
$$
where the right-hand side is holomorphic in $\Delta$. To prove the holomorphicity of ${\bf c}_j$ in $\Delta$, let us consider a polynomial (in $\lambda$)
$$
P(x_2,\ldots,x_m,\lambda):=\det(A-\lambda I)=(-\lambda)^n+({\rm tr}\,A)(-\lambda)^{n-1}+\ldots+\det A, \quad (x_2,\ldots,x_m)\in\Delta.
$$
Since its coefficients are holomorphic in $\Delta$, for any $(x_2,\ldots,x_m)\in\Delta$ the roots of $P(x_2,\ldots,x_m,\lambda)$ are close to those of
$P(0,\ldots,0,\lambda)$. Hence, there exists $j_0\in{\mathbb N}$ such that $\det(A-jI)\ne0$ for any $j>j_0$, $(x_2,\ldots,x_m)\in\Delta$. Therefore
${\bf c}_j$ is holomorphic in $\Delta$ except, maybe (if $j\leqslant j_0$), along a complex hypersurface
$$
\{P(x_2,\ldots,x_m,j)=0\}\subset\Delta.
$$
But ${\bf c}_j$ is represented by a (formal) Taylor series at the origin, thus it is holomorphic in $\Delta$ (in fact, maybe in some smaller polydisc,
however this sequence of polydiscs of holomorphicity is stabilized for $j>j_0$).

Now we can regard the (vector) formal series $\psi$ as a formal solution of an ODEs system of Fuchs type
$$
x_1\,\frac{d y}{d x_1}={\bf f}_1(x_1,z,{\bf c}_0(z)+y),
$$
depending on the parameter $z=(x_2,\ldots,x_m)\in{\mathbb C}^{m-1}$: $\psi=\sum_{j=1}^{\infty}{\bf c}_j(z)\,x_1^j$, where
the coefficients ${\bf c}_j$ are holomorphic in the common polydisc $\Delta$. As $\psi(0)=0$, the convergence of $\psi$ follows from a natural
generalization (for the multidimensional parameter $z$) of the lemma we present below.
\medskip

{\bf Lemma 1} (\cite{GS}, Ch. III, Lemme 1.1). {\it Consider an ODEs system
$$
x_1\,\frac{d y}{d x_1}={\bf f}(x_1,z,y),
$$
where $z\in\mathbb C$ is a parameter and ${\bf f}:\,({\mathbb C}\times{\mathbb C}\times{\mathbb C}^n,0)\rightarrow({\mathbb C}^n,0)$ is a germ of a
holomorphic map. Then its any formal solution $\sum_{j=0}^{\infty}{\bf a}_j(z)\,x_1^j$, where all ${\bf a}_j$ are holomorphic in some common disc
centered at the origin and ${\bf a}_0(0)=0$, is convergent.
}

\section{Complement to Theorem C. The bivariate case}

For any $f$ from ${\mathbb C}[[x]]$ (respectively, from ${\mathbb C}[[x]]^n$ or ${\rm Mat}(n,{\mathbb C}[[x]])$) we will say that
$$
{\rm ord}_{x_i}f\geqslant p\in{\mathbb Z}_+,
$$
if $f=x_i^p\,g$, with $g$ from ${\mathbb C}[[x]]$ (respectively, from ${\mathbb C}[[x]]^n$ or ${\rm Mat}(n,{\mathbb C}[[x]])$).

In this section we will consider the case of two independent variables $x_1$, $x_2$ (that is, $m=2$).

\medskip

{\bf Theorem 4.} {\it If the Jacobi matrices $\partial{\bf f}_i/\partial y$ on the formal solution $(\ref{series})$ of
$(\ref{PDE})$ satisfy
$$
{\rm ord}_{x_i}\frac{\partial{\bf f}_i}{\partial y}(x_1,x_2,\varphi)\geqslant p_i-1, \qquad i=1,\,2,
$$
then the series $\varphi$ converges near $0\in{\mathbb C}^2$.}
\medskip

We start the proof of Theorem 4 with a lemma for an ODEs system
\begin{equation}\label{ODEparam}
x_1^p\,\frac{d y}{d x_1}={\bf f}(x_1,z,y), \qquad y=(y_1,\ldots,y_n)^{\top},
\end{equation}
depending on the parameter $z\in{\mathbb C}$, where
${\bf f}:\,({\mathbb C}\times{\mathbb C}\times{\mathbb C}^n,0)\rightarrow({\mathbb C}^n,0)$
is a germ of a holomorphic map and $p\geqslant1$ is an integer.
\medskip

{\bf Lemma 2.} {\it Let
$$
\varphi=\sum_{j=0}^{\infty}{\bf a}_j(z)\,x_1^j, \qquad {\bf a}_0(0)=0,
$$
be a formal power series solution of $(\ref{ODEparam})$ with the vector coefficients ${\bf a}_j$ holomorphic in some
common disc $\Delta\subset{\mathbb C}$ of the parameter space centered at the origin. If
$$
{\rm ord}_0\frac{\partial{\bf f}}{\partial y}(x_1,z,\varphi)\geqslant p-1,
$$
then $\varphi$ converges in a neighbourhood of zero.}
\medskip

{\bf Proof.}  Let us represent $\varphi$ in the form
$$
\varphi=\sum_{j=0}^N{\bf a}_j(z)\,x_1^j+x_1^N({\bf a}_{N+1}(z)x_1+\ldots)=:\varphi_N+x_1^N\psi,
$$
with $N\geqslant p-1$. Then the formal power series $\psi$ will satisfy
\begin{equation}\label{psirelation}
x_1^p\,\frac{d\varphi_N}{d x_1}+x_1^{N+p}\frac{d\psi}{d x_1}+Nx_1^{N+p-1}\psi={\bf f}(x_1,z,\varphi_N)+
x_1^N\frac{\partial{\bf f}}{\partial y}(x_1,z,\varphi_N)\psi+x_1^{2N}O(\psi^2).
\end{equation}
Since
$$
\frac{\partial{\bf f}}{\partial y}(x_1,z,\varphi)=\frac{\partial{\bf f}}{\partial y}(x_1,z,\varphi_N)+x_1^NO(\psi),
\qquad {\rm ord}_0\psi\geqslant1,
$$
the assumption of the lemma also implies
$$
{\rm ord}_0\frac{\partial{\bf f}}{\partial y}(x_1,z,\varphi_N)\geqslant p-1,
$$
hence due to (\ref{psirelation}), the order of ${\bf f}(x_1,z,\varphi_N)-x_1^p\,d\varphi_N/d x_1$ is not
less than $N+p-1$. This means that we can devide the relation (\ref{psirelation}) by $x_1^{N+p-1}$ and obtain
the ODEs system of Fuchs type satisfied by the formal power series $\psi$:
$$
x_1\,\frac{d\psi}{d x_1}={\bf h}(x_1,z,\psi),
$$
where ${\bf h}:\,({\mathbb C}\times{\mathbb C}\times{\mathbb C}^n,0)\rightarrow({\mathbb C}^n,0)$
is a germ of a holomorphic map. Since $\psi(0)=0$, the convergence of $\psi$, as in Theorem 1, again follows from Lemma 1.
This proves Lemma 2.
\medskip

Let us show that the formal power series solution $\varphi$ can be represented in the form
$$
\varphi=\sum_{j=0}^{\infty}{\bf c}_j(x_2)\,x_1^j,
$$
where ${\bf c}_j\in{\mathbb C}\{x_2\}^n$ are holomorphic vector functions in some common disc
$\Delta\subset{\mathbb C}$ centered at the origin. As in Theorem 1, $\varphi$ is represented in such a form with
${\bf c}_j\in{\mathbb C}[[x_2]]^n$, and the main task is to prove the holomorphicity (convergence)
of ${\bf c}_j$. First one represents $\varphi$ as
$$
\varphi={\bf c}_0+\psi, \qquad \psi=\sum_{j=1}^{\infty}{\bf c}_j(x_2)\,x_1^j.
$$
Then, in view of (\ref{PDE}), ${\bf c}_0$ satisfies the relation
$$
x_2^{p_2}\,\frac{d{\bf c}_0}{d x_2}={\bf f}_2(0,x_2,{\bf c}_0),
$$
moreover, one has
$$
{\rm ord}_{x_2}\frac{\partial{\bf f}_2}{\partial y}(0,x_2,{\bf c}_0)=
{\rm ord}_{x_2}\frac{\partial{\bf f}_2}{\partial y}(x_1,x_2,\varphi)|_{x_1=0}\geqslant
{\rm ord}_{x_2}\frac{\partial{\bf f}_2}{\partial y}(x_1,x_2,\varphi)\geqslant p_2-1.
$$
Thus, by Lemma 2, the formal power series ${\bf c}_0$ converges near the origin.

Further, making the change of variable $y={\bf c}_0+u$ in the system (\ref{PDE}) we obtain
an ODEs system for the unknown $u$:
$$
x_2^{p_2}\,\frac{d u}{\partial x_2}=f_2(x_1,x_2,c_0)+\frac{\partial{\bf f}_2}{\partial y}(x_1,x_2,{\bf c}_0)u+O(u^2).
$$
This has the formal power series solution $\psi=\sum_{j=1}^{\infty}{\bf c}_j(x_2)\,x_1^j$, which implies the
relations for each coefficient ${\bf c}_j$, $j\geqslant 1$:
$$
x_2^{p_2}\,\frac{d{\bf c}_j}{d x_2}=\frac{\partial{\bf f}_2}{\partial y}(0,x_2,{\bf c}_0){\bf c}_j+
{\bf g}_j(x_2,{\bf c}_0,\ldots,{\bf c}_{j-1}).
$$
As shown above,
$$
{\rm ord}_{x_2}\frac{\partial{\bf f}_2}{\partial y}(0,x_2,{\bf c}_0)\geqslant p_2-1,
$$
hence by Lemma 2, all the ${\bf c}_j$ converge near $0\in{\mathbb C}$. Since they satisfy
the linear ODEs systems with the same homogeneous part, they are also holomorphic in $\Delta
\subset{\mathbb C}$, where $\Delta$ is a common disc of the holomorphicity for
$\partial{\bf f}_2/\partial y(0,x_2,{\bf c}_0)$, ${\bf g}_1(x_2,{\bf c}_0)$ centered at
$0\in{\mathbb C}$ ($i=2,\ldots,m$).

We have that the series
$$
\varphi=\sum_{j=0}^{\infty}{\bf c}_j(x_2)\,x_1^j
$$
with the coefficients ${\bf c}_j$ holomorphically depending on the parameter $z=x_2\in\Delta$, satisfies the ODEs
system
$$
x_1^{p_1}\,\frac{dy}{dx_1}={\bf f}_1(x_1,z,y),
$$
furthermore
$$
{\rm ord}_{x_1}\frac{\partial{\bf f}_1}{\partial y}(x_1,z,\varphi)\geqslant p_1-1.
$$
Applying Lemma 2 one proves the convergence of $\varphi$ near $0\in{\mathbb C}^2$.
\medskip

\medskip
The work of the first author is supported by the Russian Foundation for Basic Research under grant no. 18-01-00422 А
and the Program of the Presidium of the Russian Academy of Sciences №01 'Fundamental Mathematics and its Applications' under grant PRAS-18-01.




\end{document}